\documentclass[10pt,reqno]{amsart}
\usepackage{amssymb,amscd,amsbsy}
\usepackage{amsthm}
\newcommand{\z}{\zeta}

\renewcommand{\l}{\lambda}

\renewcommand{\t}{\tau}

\renewcommand{\o}{\omega}

\newcommand{\D}{\Delta}

\newcommand{\E}{{\mathcal E}}

\newcommand{\h}{{\mathcal H}}

\newcommand{\K}{{\mathcal K}}

\newcommand{\C}{{\Bbb C}}
\newcommand{\T}{{\Bbb T}}

\newcommand{\dd}{{\Bbb D}}
\newcommand{\R}{{\Bbb R}}
\newcommand{\Z}{{\Bbb Z}}

\newcommand{\bs}{\boldsymbol}

\newcommand{\bS}{{\boldsymbol S}}

\newcommand{\rf}[1]{(\ref{#1})}

\newcommand{\df}{\stackrel{\mathrm{def}}{=}}

\newcommand{\trace}{\operatorname{trace}}

\newcommand{\const}{\operatorname{const}}

\newcommand{\eeq}{\end{equation}}
\newcommand{\beq}{\begin{equation}}
\newcommand{\bay}{\begin{eqnarray}}
\newcommand{\ba}{\begin{align*}}
\newcommand{\ea}{\end{align*}}
\newcommand{\ey}{\end{eqnarray}}
\newcommand{\bey}{\begin{eqnarray*}}
\newcommand{\eey}{\end{eqnarray*}}

\newcommand{\be}{\infty}

\newcommand{\im}{\operatorname{Im}}

\newtheorem{thm}{\hspace{\parindent}Theorem}[section]

\theoremstyle{remark}

\newtheorem*{rem*}{Remark}

\newcommand\dg{\frak D}

\newcommand\OL{\operatorname{OL}}
\newcommand{\OLA}{{\rm OL}_{\rm A}}
\newcommand{\vn}{\varnothing}

\begin{document}

\newcommand{\vse}{\vspace{.2in}}

\title{A trace formula for functions of contractions and analytic operator Lipschitz functions}

\maketitle
\begin{center}
\Large
Mark Malamud$^{\rm a}$, Hagen Neidhardt$^{\rm b}$, Vladimir Peller$^{\rm c}$
\end{center}

\begin{center}
\footnotesize
{\it$^{\rm a}$Institute of Applied Mathematics and Mechanics, NAS of Ukraine,
Slavyansk, Ukraine,
and RUDN University, 6 Miklukho-Maklay St., Moscow, 117198, Russia \\
$^{\rm b}$Institut f\"ur Angewandte Analysis und Stochastik,
Mohrenstr. 39, D-10117 Berlin, Germany
\\
$^{\rm c}$Department of Mathematics, Michigan State University, East Lansing, MI 48824, USA and\\ RUDN University, 6 Miklukho-Maklay St., Moscow, 117198, Russia}
\end{center}

\newcommand{\mt}{{\mathcal T}}

\footnotesize

{\bf Abstract.} In this note we study the problem of evaluating  the trace of $f(T)-f(R)$, where $T$ and $R$ are contractions on Hilbert space with trace class difference, i.e., $T-R\in\bS_1$ and $f$ is a function analytic in the unit disk 
$\dd$. It is well known that if $f$ is an operator Lipschitz function analytic in $\dd$, then $f(T)-f(R)\in\bS_1$. The main result of the note says that there exists a function $\bs{\xi}$ (a spectral shift function) on the unit circle $\T$ of class $L^1(\T)$ such that the following trace formula holds: $\trace(f(T)-f(R))=\int_\T f'(\z)\bs{\xi}(\z)\,d\z$, whenever $T$ and $R$ are contractions with $T-R\in\bS_1$ and $f$ is an operator Lipschitz function analytic in $\dd$.

\

\begin{center}
{\bf\large Une formule de trace pour les fonctions de contractions et
les fonctions analytique lipschitziennes op\'eratorilelles}
\end{center}

\medskip

{\bf R\'esum\'e.} Nous consid\'erons dans cette note le probl\`eme de trouver le trace de $f(T)-f(R)$ o\`u $T$ et $R$ sont des contractions dans un espace hilbertien et $f$ est une fonction analytique dans le disque unit\'e $\dd$. Il est bien connu que si $f$ est une fonction analytique dans $\dd$ qui est lipschitzienne op\'eratorilelle est la diff\'erence $T-R$ est de classe trace, c'est-\`a-dire
$T-R\in\bS_1$, alors $f(T)-f(R)\in\bS_1$. Le r\'esultat principal de cette note \'etablit qu'il existe une fonction $\bs{\xi}$ (une fonction de d\'ecalage spectral) sur le cercle unit\'e $\T$ dans l'espace $L^1(\T)$ pour laquelle la formule de trace suivante est vrai:
$\trace(f(T)-f(R))=\int_\T f'(\z)\bs{\xi}(\z)\,d\z$ pour n'importe quelle
fonction $f$ lipschitzienne op\'eratorilelle et analytique dans $\dd$.
\normalsize

\medskip

\begin{center}
{\bf\large Version fran\c caise abr\'eg\'ee}
\end{center}

\medskip

La fonction de d\'ecalage spectral pour des couple d'op\'erateurs auto-adjoints
\'etait introduit par I.M. Lifshits dans \cite{L}. M.G. Krein consid\'erait dans \cite{Kr} le cas le plus g\'en\'eral. Soient $A$ et $B$ des op\'erateurs auto-adjoints (pas n\'ecessairement born\'es) dont la diff\'erence $A-B$ est de classe trace, c'est-\`a-dire $A-B\in\bS_1$. Il \'etait d\'emontr\'e dans \cite{Kr} qu'il existe une fonction $\bs{\xi}=\bs{\xi}_{A,B}$ r\'eelle dans $L^1(\R)$ (qui depend de $A$ et $B$) pour laquelle la formule de trace suivante est vrai:
\bay
\label{foslsso}
\trace\big(f(A)-f(B)\big)=\int_\R f'(t)\bs{\xi}_{A,B}(t)\,dt
\ey
pour chaque fonction $f$ diff\'erentiable sur $\R$ telle que la d\'eriv\'ee $f'$ de $f$ est la transform\'ee de Fourier d'une mesure complexe borelienne sur $\R$. 
La fonction $\bs{\xi}$ s'appelle la {\it fonction de d\'ecalage spectral pour la couple} $(A,B)$.
M.G. Krein a pos\'e dans \cite{Kr} le probl\`eme de d\'ecrire la classe de fonctions $f$ pour lesquelles la formule de trace ci-dessus est vrai pour toutes les couples 
d'op\'erateurs auto-adjoints $(A,B)$ telles que $A-B\in\bS_1$.

Le probl\`eme de Krein \'etait r\'esolu r\'ecemment dans \cite{Pe4}: la classe de fonctions ci-dessus co\"incide avec la classe de fonctions lipschitziennes op\'eratorilelles sur $\R$. Rappelons qu'une fonction $f$ continue sur $\R$ s'appelle {\it lipschitzienne op\'eratorielle} si on a
\bay
\label{opLip}
\|f(A)-f(B)\|\le\const\|A-B\|
\ey
pour tous les op\'erateurs auto-adjoints $A$ et $B$.

Dans le travail \cite{Kr2} M.G. Krein a introduit la fonction de d\'ecalage spectral pour les couple d'op\'erateurs unitaires dont la diff\'erence est de classe trace. Il a d\'emontr\'e que pour chaque couple $(U,V)$ d'op\'erateurs unitaires pour lesquels $U-V\in\bS_1$ il existe une fonction $\bs{\xi}_{U,V}$ dans l'espace $L^1(\T)$ (qui s'appelle une fonction de d\'ecalage spectral pour la couple $(U,V)$) telle que
\bay
\label{fosuo}
\trace\big(f(U)-f(V)\big)=\int_\T f'(\z)\bs{\xi}_{U,V}(\z)\,d\z
\ey
pour chaque fonction $f$ diff\'erentiable dont la d\'eriv\'ee a une s\'eries de Fourier absolument convergente.

Le probl\`eme de d\'ecrire la classe maximale de fonctions $f$ pour lesquelles la formule \rf{fosuo} s'applique pour toutes les couples $(U,V)$ d'op\'erateurs unitaires avec $U-V\in\bS_1$ \'etait r\'esolu r\'ecemment dans \cite{AP2}.
Notamment, il \'etait d\'emontr\'e dans \cite{AP2} que la classe dont il s'agit co\"incide avec la classe de fonctions lipschitziennes op\'eratorielles sur le cercle $\T$.

Dans cette note nous consid\'erons le cas de fonctions des contractions sur l'espace hilbertien. Rappelons q'on dit qu'un op\'erateur $T$ sur l'espace hilbertien s'appelle une {\it contaction} si $\|T\|\le1$.

Le r\'esultat principal de cette note est le th\'eor\`eme suivant:

\medskip

{\bf Th\'eor\`eme.} {\it Pour chaque couple $(T,R)$ de contractions sur l'espace hilbertien il existe une fonction $\bs{\xi}=\bs{\xi}_{T,R}$ de l'espace $L^1(\T)$ (une fonction de d\'ecalage spectral pour $T$ et $R$) pour laquelle la formule de trace suivante
\bay
\label{fosszha}
\trace\big(f(T)-f(R)\big)=\int_\T f'(\z)\bs{\xi}(\z)\,d\z
\ey
s'applique pour toutes les fonctions $f$ lipschitziennes op\'eratorielles et analytique dans $\dd$.}

\medskip

Remarquons que la classe de fonctions lipschitziennes op\'eratorielles et analytique dans $\dd$ est la classe maximale de fonctions pour lesquelles la formule \rf{fosszha} est vrai pour tous le contractions $T$ et $R$ dont la diff\'erence est de classe trace.

\medskip

\begin{center}
------------------------------
\end{center}

\section{\bf Introduction}

\medskip

The notion of spectral shift function was introduced by physicist 
I.M. Lifshits in \cite{L}. It was M.G. Krein who generalized in \cite{Kr} this notion to a most general situation. Namely, if $A$ and $B$ are (not necessarily bounded) self-adjoint operators on Hilbert space with trace class difference (i.e., $A-B\in\bS_1$), then it was shown in \cite{Kr} that there exists a unique real function $\bs{\xi}=\bs{\xi}_{A,B}$ in $L^1(\R)$, the {\it spectral shift function for the pair $(A,B)$}, such that trace formula \rf{foslsso} holds
for all functions $f$ that are differentiable on $\R$ and whose derivative $f'$ is the Fourier transform of a complex Borel measure. 

Krein observed in \cite{Kr} that the right-hand side of \rf{foslsso} makes sense for arbitrary Lipschitz functions $f$ and he posed the problem to describe the maximal class of functions $f$, for which trace formula \rf{foslsso} holds for an arbitrary pair $(A,B)$ of self-adjoint operators with $A-B\in\bS_1$. 

It was Farforovskaya who proved in \cite{F} that there exist self-adjoint operators $A$ and $B$ with $A-B\in\bS_1$ and a Lipschitz function $f$ on $\R$ such that $f(A)-f(B)\not\in\bS_1$. Thus, trace formula \rf{foslsso} cannot be generalized to the class of all Lipschitz functions $f$. 
In \cite{Pe1} and \cite{Pe2} it was shown that trace formula \rf{foslsso} holds for all functions $f$ in the (homogeneous) Besov class $B_{\be,1}^1(\R)$. 

Krein's problem was completely solved recently in \cite{Pe4}. It was shown in \cite{Pe4} that the maximal class of functions $f$, for which \rf{foslsso} holds whenever $A$ and $B$ are (not necessarily bounded) self-adjoint operators with trace class difference coincides with the class of operator Lipschitz functions $f$ on $\R$. Recall that $f$ is called an {\it operator Lipschitz function} if
inequality \rf{opLip} holds for arbitrary self-adjoint operators $A$ and $B$.
We refer the reader to \cite{AP1} for detailed information on operator Lipschitz functions.

Later M.G. Krein introduced in \cite{Kr2} the notion of spectral shift function for pairs of unitary operators with trace class difference. He proved that for a pair $(U,V)$ of unitary operators with $U-V\in\bS_1$, there exists a function $\bs{\xi}=\bs{\xi}_{U,V}$ in $L^1(\T)$ (a {\it spectral shift function for the pair $(U,V)$}) such that trace formula \rf{fosuo} holds for an arbitrary differentiable function $f$ on the unit circle $\T$ whose derivative has absolutely convergent Fourier series. Note that $\bs\xi$ is unique modulo a constant additive;
it can be normalized by the condition $\int_\T\bs{\xi}(\z)\,|d\z|=0$.

An analog of the result of \cite{Pe4} was obtained in \cite{AP2}. It was proved in 
\cite{AP2} that the maximal class of functions $f$, for which trace formula \rf{fosuo} holds for arbitrary unitary operators $U$ and $V$ with trace class difference coincides with the class of operator Lipschitz functions on the unit circle; this class can be defined by analogy with operator Lipschitz functions on $\R$.
Note that the method used in \cite{Pe4}  does not work in the case of unitary operators.
We denote the class of operator Lipschitz functions on $\T$ by $\OL_\T$.

In this note we consider the case of functions of contractions. Recall that an operator $T$ on Hilbert space is called a {\it contraction} if $\|T\|\le1$.
For a contraction $T$, the Sz.-Nagy--Foia\c s functional calculus associates with each function $f$ in the disk-algebra $C_{\rm A}$ the operator $f(T)$. The functional calculus $f\mapsto f(T)$ is linear and multiplicative and 
$
\|f(T)\|\le\max\{|f(\z)|:~\z\in\C,~|\z|\le1\}
$
(von Neumann's inequality). As usual, $C_{\rm A}$ stands for the space of functions analytic in the unit disk $\dd$ and continuous in the closed unit disk. The purpose of this note is to obtain analogs of the above mentioned results of \cite{Kr}, \cite{Kr2}, \cite{Pe4} and \cite{AP2} for functions of contractions.

We are going to prove the existence of a spectral shift function for pairs $(T_0,T_1)$ of contractions with trace class difference. This is an integrable function $\bs{\xi}$ on the unit circle $\T$ such that 
\bay
\label{tfszha}
\trace\big(f(T_1)-f(T_0)\big)=\int_\T f'(\z)\bs{\xi}(\z)\,d\z
\ey
for all analytic polynomials $f$. Such a function $\bs{\xi}$ is called a {\it spectral shift function} for the pair $(T_0,T_1)$. Is is unique up to an additive in the Hardy class $H^1$. In other words, if $\bs{\xi}$ is a spectral shift function for $(T_0,T_1)$, then all spectral shift functions for the pair $(T_0,T_1)$ are given by $\{\bs{\xi}+h:~h\in H^1\}$.

The second principal result of this note is that the maximal class of functions $f$ in $C_{\rm A}$, for which formula \rf{tfszha} holds for all such pairs $(T_0,T_1)$ coincides with the class of operator Lipschitz functions analytic in $\dd$. We say that a function $f$ analytic in $\dd$ is called 
{\it operator Lipschitz} if 
$$
\|f(T)-f(R)\|\le\const\|T-R\|
$$
for contractions $T$ and $R$. 
We denote the class of operator Lipschitz functions analytic in $\dd$ by $\OLA$.
It is well known that if $f\in\OLA$, then $f\in C_{\rm A}$ and
$\OLA=\OL_\T\bigcap C_{\rm A}$ (see \cite{KS} and \cite{AP1}).

It turns out that as in the case of functions of self-adjoint operators and functions of unitary operators, the maximal class of functions, for which trace formula \rf{tfszha} holds for all pairs of contractions $(T_0,T_1)$ with trace class difference coincides with the class $\OLA$.

To obtain the results described above, we combine two approaches. The first approach is based on double operator integrals with respect to semi-spectral measures. It leads to a trace formula 
$\trace\big(f(T)-f(R)\big)=\int_\T f'(\z)\,d\nu(\z)$ for a Borel measure $\nu$ on $\T$.

The second approach is based on an improvement of a trace formula obtained in \cite{MN} for functions of dissipative operators.

\

\section{\bf Double operator integrals and a trace formula for arbitrary functions in $\bs{\OLA}$}

\medskip

Double operator integrals
$$
\iint\Phi(x,y)\,dE_1(x)Q\,dE_2(y)
$$
were introduced by Birman and Solomyak in \cite{BS}. Here $\Phi$ is a bounded measurable function, $E_1$ and $E_2$ are spectral measures on Hilbert space and
$Q$ is a bounded linear operator. Such double operator integrals are defined for arbitrary bounded measurable functions $\Phi$ if $Q$ is a Hilbert--Schmidt operator. If $Q$ is an arbitrary bounded operator, then for the double operator integral to make sense, $\Phi$ has to be a Schur multiplier with respect to $E_1$ and $E_2$, (see \cite{Pe1} and \cite{AP1}).

In this note we deal with double operator integrals with respect to {\it semi-spectral measures}
$$
\iint\Phi(x,y)\,d\E_1(x)Q\,d\E_2(y).
$$
Such double operator integrals were introduced in \cite{Pe3} (see also \cite{Pe+}).
We refer the reader to recent paper \cite{AP1} for detailed information about double operator integrals.

If $T$ is a contraction on a Hilbert space $\h$, it has a {\it minimal unitary dilation} $U$, i.e., $U$ is a unitary operator on a Hilbert space $\K$, $\K\supset\h$, $T^n=P_\h U^n\big|\h$ for $n\ge0$ and $\K$ is the closed linear span of
$U^n\h$, $n\in\Z$ (see \cite{SNF}). Here $P_\h$ is the orthogonal projection onto $\h$.
The {\it semi-spectral measure $\E_T$ of} $T$ is defined by
$$
\E_T(\D)\df P_\h E_U(\D)\big|\h,
$$
where $E_U$ is the spectral measure of $U$ and $\D$ is a Borel subset of $\T$. It is well known that 
$
T^n=\int_\T\z^n\,d\E_T(\z),~ n\ge0.
$

If $f\in\OLA$, then the divided difference
$\dg f$, 
$$
(\dg f)(\z,\t)\df\big(f(\z)-f(\t)\big)(\z-\t)^{-1},\quad\z,~\t\in\T,
$$
is a Schur multiplier with respect to arbitrary Borel (semi-)spectral measures on $\T$
and 
$$
f(T_1)-f(T_0)=\iint_{\T\times\T}(\dg f)(\z,\t)\,d\E_{T_1}(\z)(T_1-T_0)\,d\E_{T_0}(\t)
$$
for an arbitrary pair of contractions $(T_0,T_1)$ 
with trace class difference, see \cite{AP1}.


\begin{thm}
Let $f\in\OLA$ and let $T_0$ and $T_1$ be contractions on Hilbert space and
$T_t=T+t(R-T)$, $0\le t\le1$. Then
\bay
\label{proizvo}
\lim_{s\to0}\frac1s\big(f(T_{t+s})-f(T_t)\big)
=\iint_{\T\times\T}\big(\dg f)(\z,\t)\,d\E_t(\z)(T_1-T_0)\,d\E_t(\t)
\ey
in the strong operator topology,
where $\E_t$ is the semi-spectral measure of $T_t$.
\end{thm}

It can be shown that if $T_1-T_0\in\bS_1$, then 
$$
f(T_1)-f(T_0)=\int_0^1Q_t\,dt,
$$
where $Q_t$ is the right-hand side of \rf{proizvo}, and  $Q_t\in\bS_1$ for every $t\in[0,1]$. The integral can be understood in the sense of Bochner in the space $\bS_1$. It can be shown that $\trace Q_t=\int_\T f'(\z)\,d\nu_t(\z)$, where $\nu_t$ is defined by $\nu_t(\D)\df\trace\big((T-R)\E_t(\D)\big)$.
We can define now the Borel measure $\nu$ on $\T$ by
\bay
\label{oprnu}
\nu\df\int_0^1\nu_t\,dt,
\ey
which can be understood as the integral of the vector-function $t\mapsto\nu_t$ that is continuous in the weak-star topology in the space of complex Borel measures on $\T$.


\begin{thm}
\label{pomere}
Let $T_0$ and $T_1$ be contractions on Hilbert space such that $T_1-T_0\in\bS_1$.
Then 
\bay
\label{formsledlyaszha}
\trace\big(f(T_1)-f(T_0)\big)=\int\limits_\T f'(\z)\,d\nu(\z)
\ey
for every $f$ in $\OLA$, where $\nu$ is the Borel measure defined by {\em\rf{oprnu}}.
\end{thm}

\

\section{\bf A spectral shift function for a pair of contractions with trace class difference}

\medskip

In this section we obtain the existence of a spectral shift function for pairs of contractions with trace class difference. 

\begin{thm}
\label{fspsdszh}
Let $T_0$ and $T_1$ be contractions on Hilbert space with trace class difference.
Then there exists a complex function $\bs{\xi}$ in $L^1(\T)$ such that
for an arbitrary analytic polynomial $f$,
\bay
\label{T-l}
\trace\big(f(T_1)-f(T_0)\big)=
\int_\T f'(\z)\bs{\xi}(\z)\,d\z.
\ey
Moreover, if $T_0$ is a unitary operator, we can find such a function $\bs{\xi}$ that also satisfies the requirement $\im\bs{\xi}\le0$. On the other hand, if $T_1$ is a unitary operator,  we can add the requirement $\im\bs{\xi}\ge0$.
\end{thm}

%

\medskip

{\bf Remark.}
{\it It is not true in general that for a pair of contractions with trace class difference, there exists a real spectral shift function}. However, this is true under certain assumptions. In particular, if $\bs\xi$ is a spectral shift function and $\bs{\xi}\log(1+|\bs{\xi}|)\in L^1(\T)$, then we can find a real spectral shift function for the same pair of contractions.
The same conclusion holds if $\xi$ is a spectral shift function that belongs to
the weighted space $L^p(\T,w)$, where $1<p<\be$ and $w$ satisfies the Muckenhoupt condition $(A_p)$.

\medskip

To prove Theorem \ref{fspsdszh}, we can improve Theorem 3.14 of \cite{MN} and deduce Theorem \ref{fspsdszh} from that improvement with the help of Cayley transform.
On the other hand, Theorem \ref{fspsdszh} allows us to obtain a further improvement of Theorem 3.14 of \cite{MN} and obtain the following result:


\begin{thm}
\label{uluch}
Let $L_0$ and $L_1$ be maximal dissipative operators such that
\bay
\label{yadrez}
(L_1+{\rm i}I)^{-1}-(L_0+{\rm i}I)^{-1}\in\bS_1.
\ey
Then there exists a complex measurable function $\bs{\o}$  (a spectral shift function for $(L_0,L_1)$) such that
\bay
\label{vesovL1}
\int_\R|\bs{\o}(t)|(1+t^2)^{-1}\,dt<\be,
\ey
for which the following trace formula holds:
\bay
\label{fsleddis}
\trace\big((L_1-\l I)^{-1}-(L_0-\l I)^{-1}\big)
=-\int_\R\bs{\o}(t)(t-\l)^{-2}\,dt,\quad\im\l<0.
\ey
Moreover, if $L_0$ is self-adjoint, there exist a function $\bs{\o}$ satisfying
{\em\rf{vesovL1}} and {\em\rf{fsleddis}} such that 
$\im\bs{\o}\ge0$ on $\R$, while if $L_1$ is self-adjoint, there exist a function $\bs{\o}$ satisfying
{\em\rf{vesovL1}} and {\em\rf{fsleddis}} and that 
$\im\bs{\o}\le0$ on $\R$.
\end{thm}

Recall that a closed densely defined operator $L$ is called {\it dissipative} if $\im(Lx,x)\ge0$ for every $x$ in its domain. It is called a {\it maximal dissipative operator} if it does not have a proper dissipative extension.


\medskip

{\bf Remark.} In the case when $L_0-L_1\in\bS_1$, Theorem \ref{uluch} can be specified. Namely, it was shown in \cite{MN} (Theorem 4.11) that a spectral shift function $\bs{\o}$ can be chosen in $L^1(\R)$.

\medskip

Note also that Theorem \ref{fspsdszh} improves earlier results in \cite{AN} and \cite{R2},
while Theorem \ref{uluch} improves Theorem 3.14 of \cite{MN} (the latter imposes the additional assumption $\rho(L_0)\cap\C_+\neq\vn$)
and also improves and complements earlier results in \cite{R1} and \cite{Kr3} (see \cite{MN} for details).

\

\section{\bf The main result}

\medskip

Now we are able to state the main result of this note. 

\begin{thm}
\label{osnteor}
Let $T_0$ and $T_1$ be contractions satisfying $T_1-T_0\in\bS_1$ and let $\bs{\xi}$ be a spectral shift function for $(T_0,T_1)$. Then for every $f\in\OLA$
the following trace formula holds
\bay
\label{zavfo}
\trace\big(f(T_1)-f(T_0)\big)=\int_\T f'(\z)\bs{\xi}(\z)\,d\z.
\ey
\end{thm}

Indeed, by Theorem 
\ref{fspsdszh}, formula \rf{zavfo} holds for analytic polynomials  $f$. Combining this fact with formula \rf{formsledlyaszha}, we see that the measure $\nu$ is absolutely continuous with respect to normalized Lebesgue measure and differs from the measure $\bs{\xi}dz$ by an absolutely continuous measure with Radon--Nikodym density in $H^1$.

\medskip

{\bf Remark.} It is easy to see that the condition that $f$ has to be operator Lipschitz is not only sufficient for formula \rf{zavfo} to hold for arbitrary pairs of contractions $(T_0,T_1)$ with trace class difference, but also necessary. Indeed, it is well known (see \cite{AP1}) that if $f$ is not operator Lipschitz, then there exist unitary operators $U$ and $V$ such that $U-V\in\bS_1$, but $f(U)-f(V)\not\in\bS_1$.

\medskip

By applying Cayley transform, we can deduce now from Theorem \ref{osnteor} the following analog of it for dissipative operators.

\begin{thm}
Let $L_0$ and $L_1$ be maximal dissipative operators satisfying
{\em\rf{yadrez}}. Suppose that $f$ is a function analytic in the upper half-plane snd such that the function 
$$
\z\mapsto f\big(({\rm i}-\z)({\rm i}+\z)^{-1}\big),\quad\z\in\dd,
$$
belongs to $\OLA$. Then $f(L_1)-f(L_0)\in\bS_1$ and
\bay
\label{sledrraya}
\trace\big(f(L_1)-f(L_0)\big)=\int_\R f'(t)\bs{\o}(t)\,dt,
\ey
where $\bs\o$ is a spectral shift function for the pair $(L_0,L_1)$.
\end{thm}

\medskip

{\bf Remark.} In the case when $L_1-L_0\in\bS_1$ and $\bs\o\in L^1(\R)$, it can be shown that formula \rf{sledrraya} holds for all operator Lipschitz functions in the upper half-plane (see \cite{AP1} for a discussion of the class of such functions).

\medskip

The research of the first author is partially supported by by the Ministry of Education and Science of the Russian Federation
      (the Agreement number N$^{\underline\circ}$ 02.à03.21.0008), the research of the third author is partially supported by NSF grant DMS 1300924 and by the Ministry of Education and Science of the Russian Federation (the Agreement number N$^{\underline\circ}$ 02.à03.21.0008).

\

\end{document}